\documentclass[12pt]{amsart}
\usepackage{amsmath}
\usepackage{amsthm}
\usepackage{amscd}

\newcommand{\Q}{{\mathbb Q}}
\newcommand{\R}{{\mathbb R}}
\newcommand{\Z}{{\mathbb Z}}
\newcommand{\C}{{\mathbb C}}

\newcommand{\So}{{\mathbb S}}

\renewcommand{\P}{{\mathbb P}}
\newcommand{\N}{{\mathbb N}}
\DeclareMathOperator*{\Injlim}{\varinjlim}
\DeclareMathOperator*{\Projlim}{\varprojlim}

\newtheorem{theorem}{Theorem}
\newtheorem{lemma}{Lemma}
\newtheorem{corollary}{Corollary}
\newtheorem{proposition}[theorem]{Proposition}
\newtheorem*{Labesse}{Theorem (Labesse, [\cite{Lab}, Theorem 7.1])}
\theoremstyle{definition}
\newtheorem*{acknowledgement}{Acknowledgement}

\newtheorem{definition}{Definition}

\numberwithin{question}{section}
\theoremstyle{remark}
\newtheorem{remark}{Remark}

\begin{document}

\title[Cohomolgy of Weil groups]{On the vanishing of the measurable
 Schur\\  cohomology
 groups of Weil groups}

\author{C.~S.~Rajan}

\address{Tata Institute of Fundamental 
Research, Homi Bhabha Road, Bombay - 400 005, INDIA.}
\email{rajan@math.tifr.res.in}


\begin{abstract}
We generalize a theorem of Tate and show that the second cohomology of
the Weil group of a global or local field with coefficients in $\C^*$
(or more generally, with coefficients in the complex points of tori 
over $\C$) vanish, where the cohomology groups are defined using measurable 
cochains in the sense of Moore. We recover a
theorem of Labesse proving that admissible
homomorphisms of a  Weil group to the Langlands dual of a reductive group 
can be lifted to an extension of the Langlands dual group  by a tori. 

\end{abstract}

\maketitle

\section{Introduction}
A classical theorem of Tate  asserts  the vanishing of the `Schur
cohomology group (see \cite{Se1}),  
\[ H^2_{gc}(G_K, \C^*)=\{1\},\]
where $G_K$ is the absolute Galois group of  a global or local  field
$K$ and $G_K$ acts trivially on $\C^*$. The
cohomology groups are those that occur in Galois cohomology, 
constructed  using continuous cochains with
respect to the profinite topology on $G_K$ and the discrete topology
on $\C^*$. The 
obstruction to lifting a  projective representation of a finite or a
profinite group $G$  to a linear
representation of $G$ lies in the Schur cohomology group
$H^2(G,\C^*)$. Consequently, Tate's theorem asserts that any projective
representation $\tilde{\rho}:G_K\to PGL_n(\C)$ of $G_K$, 
 has a lift to a continuous, linear representation of
${\rho}:G_K\to GL_n(\C)$, satisfying $\pi\circ \rho=\tilde{\rho}$,
where  $\pi: GL_n(\C)\to PGL_n(\C)$ denotes the natural
projection.

A natural question that arises  is to establish analogues of Tate's
theorem and the lifting theorem for projective representations 
for more general `Galois type' groups. For example, instead of
$G_K$, we can consider the Weil group $W_K$ associated to a global or local
field $K$ (or more generally, the complex points of the 
conjectural Langlands group, the
continuous representations of which parametrize automorphic representations).  
These are topological groups which are no longer profinite. Since our primary motivation is the application to lifting,
we need to work with a  cohomology theory where  the second
cohomology $H^2(G,A)$ parametrizes the topological   extensions of $G$ by
$A$.  Such a
cohomology theory has been defined and developed by Moore using measurable
cochains (\cite{M1,M2,M3,M4}). 
 We show that the  analogue of Tate's theorem holds in this
framework:  
\begin{theorem}\label{gentate} Let $K$ be a local or global field, and
let  $W_K$ denote the Weil group of $K$.  Then the  Moore cohomology group  
$H^2(W_K, \C^*)$ vanishes, where the Weil group $W_K$ acts trivially
on $\C^*$, and $\C^*$ is equipped with the usual Euclidean topology. 
\end{theorem}
More generally, let $T$ be an algebraic torus over $\C$ with a
continuous action  $G_K\to \text{Aut}(T)$, where $\text{Aut}(T)$ is
the group of algebraic automorphisms of $T$. 
We let $W_K$ act on $T$ via the canonical morphism  to $G_K$, and
consider $T(\C)$ equipped with the usual Euclidean topology as a $W_K$
module. Our main
theorem is, 
\begin{theorem}\label{maintheorem}
\[H^2(W_K, T(\C))=(0).\]
\end{theorem}

One of the intriguing aspects of these vanishing theorems, is the
comparison with the corresponding Galois cohomology groups. 
Over number fields it is not true in general that $H^2_{gc}(G_K, T(\C))$
vanishes. 
 Thus the problems arising from the presence
of the (real) archimedean  places towards the Galois cohomology of complex tori
disappear upon going over to the Weil group, and the Weil group
of a number  field
behaves more like the fundamental group of a curve than the Galois
group. However, as pointed out by S. Lichtenbaum, it is not true that
the higher cohomology groups
$H^i(W_K, T(\C))$ vanish for $i\geq 3$.

Apart from the intrinsic interest, the principal motivation 
for proving the vanishing  theorem lies in the
application to lifting admissible morphisms of Weil group to the dual
Langlands group of a reductive group defined over $K$ (see Theorem
\ref{labesse}).  These 
applications to lifting have  been obtained by 
Langlands if $K$ is archimedean \cite[Lemma 2.10]{Lan3}, by Henniart over
non-archimedean local field and function fields \cite{He}, 
and by  Labesse over number fields 
\cite{Lab} (see Theorem \ref{labesse}).  For us,  the
motivation to prove such a lifting result,  is 
in the application to characterising image of solvable 
base change \cite{R} for $GL(2)$. 

The results of this paper, as pointed out by S. Lichtenbaum and the
referee,  are of relevance  to the study of the values of the 
Hasse-Weil zeta functions at integers and the concept of
Weil-\'{e}tale topology introduced by Lichtenbaum in this context.

\begin{acknowledgement}
After obtaining a proof of Theorem \ref{gentate}, 
D. Prasad brought to our attention
the results proved by Labesse \cite{Lab}, and I sincerely thank him for
this and for useful discussions. I am also indebted to S. Lichtenbaum
and the referee for pointing out the applications of the results of
this paper to Lichtenbaum's work, and to S. Lichtenbaum for showing me
that my earlier conjecture on the vanishing of the higher cohomology
groups in question is false. 
\end{acknowledgement}

\section{Weil groups}
We briefly recall the basic properties of Weil groups. The
construction and study of the basic 
properties of Weil groups has been carried out in
great detail in the notes of Artin and Tate \cite{AT}, and we refer to
this book for further details. 

Let $K$ denote a local or global field, and let 
\[
\begin{split}
C_K &= J_K/K^* \quad \text{if K is global}\\
  & = K^* \quad \text{ K local}
\end{split}
\]
The motivation for the Weil group $W_K$ 
is to obtain non-abelian generalizations of the abelian reciprocity
law. We require a completion of the diagram, 
\[
\begin{CD} 
W_K @>{\phi}>> G_K \\
@VVV  @VVV\\
C_K @>{\psi_K}>> G_K^{ab}
\end{CD}
\]
where we require that $W_K^{ab}\simeq C_K$, and 
$\psi_K:C_K\to G_K^{ab}$ is the reciprocity map of class field
theory. Here for a topological  group $G$, $G^{ab}$ denotes the
quotient of $G$ by the closure of the commutator group of $G$. 
A Weil group for $K$ is a triple $(W_K, \phi, \{r_L\})$
consisting of a topological group $W_K$, together with a continuous
homomorphism 
\[\phi: W_K\to G_K,\]
with dense image. For each finite extension $L$ of $K$,  let
$W_L=\phi^{-1}(G_L)$. We are given
isomorphisms 
\[ r_L: C_L\to W_L^{ab},\]
which composed with the map from $W_L^{ab}\to G_L^{ab}$ induced by
$\phi$, gives the reciprocity map $\psi_L$. To be a Weil group, it is
further required that these structures be compatible with changing the
field, with respect to inclusions, norm maps and conjugating by
elements of the Galois group, 
in a manner similar to the behavior of these maps with respect to  the
reciprocity map in class field
theory. In particular, this implies that $W_L$ is a Weil group
associated to $L$. 

Let $K$ now denote a non-archimedean local field or a function
field of a curve over a finite field $k$. We refer to these cases as
the $\Z$-cases.  In the first case, let $N_K$
denote the inertia group, and for function fields let $N_K$
denote the `geometric
Galois group', i.e., the absolute Galois group $G(\bar{K}/K\bar{k})$
 of a separable closure
$\bar{K}$ of $K$
over the compositum $K\bar{k}$, where $\bar{k}$ is an algebraic
closure of $k$. It is known in these cases that the reciprocity map
$\psi_K$ is injective, and that the image of $C_K$ consists of those
elements of $G(\bar{K}/K\bar{k})^{ab}$ which restrict to an integral power
of the Frobenius in $G(\bar{k}/k)$. The Weil group
admits a similar description in these cases. We have that $W_K$ is a
dense subgroup of $G_K$, and consists of those
elements of $G(\bar{K}/K\bar{k})^{ab}$ which restrict to an integral power
of the Frobenius in $G(\bar{k}/k)$. There is an exact sequence,
\begin{equation}\label{Z-case}
 1\to N_K\to W_K\to \Z \to 1.
\end{equation}
The topology on $W_K$ comes from the product topology on $N_K\times
\Z$, where $\Z$ is discrete and $N_K$ is a profinite group. 

The structure of the Weil group over archimedean fields is known. The
Weil group of $\C$ is $\C^*$, and the Weil group of $\R$ is an
extension of $\Z/2\Z$ by $\C^*$, where the non-trivial element of
$\Z/2\Z$ acts by conjugation on $\C^*$.  The $2$-cocycle $a(\tau_1, \tau_2)$ 
defining the
extension is given by 
\begin{equation} \label{extensionclass}
a(\tau_1, \tau_2)=
\begin{cases} 1 & \text{if either $\tau_1$ or
$\tau_2$ is $1$},\\ $-1$ & \text{if $\tau_1=\tau_2=-1$}.
\end{cases}
\end{equation}

For the Weil groups over number fields only a construction is
known in terms of the fundamental class, which we now recall. 
Let 
\[ H^2(G_K, C_{\bar{K}})=\Injlim_{\{L: K\subset
L\subset \bar{K}\}} H^2(G(L/K), C_L),\]
the limit being taken over the finite Galois extensions $L$ of $K$
contained inside a separable algebraic closure $\bar{K}$ of $K$. 

One of the main theorems of class field thery is that there exists a
`canonical isomorphism'
\[ {\rm inv}_K:  H^2(G_K, C_{\bar{K}})\to \Q/\Z,\]
compatible with localization and base change. If
$L/K$ is a finite Galois extension of degree $n$, we have an inclusion
\[ H^2(G(L/K), C_L)\to H^2(G_K, C_{\bar{K}})\]
and the image is the cyclic group $\frac{1}{n}\Z/\Z$ of order $n$. 
The fundamental class is defined as the element $\alpha_{L/K}\in
H^2(G(L/K), C_L)$ such that 
\[ {\rm inv}_K (\alpha_{L/K})=\frac{1}{n}.\]
Let $W_{L/K}$ be the  locally compact topological group 
defined by the class $\alpha_{L/K} \in
H^2(G(L/K), C_L)$:
\[ 1\to C_L\to W_{L/K}\to G(L/K)\to 1.\]
Given finite Galois extensions $K\subset L\subset M$, we obtain
isomorphisms of Weil  groups,
\begin{equation*}
\begin{split}
W_{M/L} & \simeq \phi^{-1}(G_{M/L})\\
{\rm and}~~~W_{L/K} & \simeq W_{M/K}/W_{M/L}^c,
\end{split}
\end{equation*}
where the notation is that for any subgroup $H$ of a topological group
$G$, $H^c$ denotes the closure of the commutator subgroup of
$H$. It follows from the axioms defining the Weil group that the
induced map $C_M\to C_L$ is given by the norm map. Further there is a
commutative diagram,
\[
\begin{CD}
1 \to C_L  @>>> W_{L/K} @>>> G(L/K) \to  1 \\
 @VVV  @V{\phi_{L/K}}VV @V{=}VV\\
1\to G(L^{ab}/L) @>>> G(L^{ab}/K) @>>> G(L/K) \to 1
\end{CD}
\]
Choosing a transitive system of Weil morphisms, the Weil group
$W_K$ is obtained as a  projective limit,
\[ W_K=\Projlim_{\{L: K\subset
L\subset \bar{K}\}} W_{L/K}.\]
The morphism $\phi$ is obtained as the projective limit of the
morphisms $\phi_{L/K}$, 
\[ \phi=\Projlim_{\substack{\{L: K\subset
L\subset \bar{K}\}}} \phi_{L/K}:\Projlim_{\{L:K\subset
L\subset \bar{K}\}} W_{L/K}\to \Projlim_{\{L: K\subset
L\subset \bar{K}\}} G(L^{ab}/K)\simeq G_K.\]
Let $N_K$ denote the kernel of $\phi$.   
Let $D_L$ denote the connected component of
identity of $C_L$. We have an exact sequence given by the reciprocity
map,
\begin{equation}\label{reciprocity}
 1\to D_L\to C_L\to G(L^{ab}/L)\to 1.
\end{equation}
The structure of $D_L$ is known. If $L$ has $r_1(L)$ real places and
$2r_2(L)$ complex places, then 
\begin{equation}\label{DL}
D_L\simeq \R\times (S^1)^{r_2(L)}\times \So^{r_1(L)+r_2(L)-1},
\end{equation}
where $S^1$ denotes  the usual circle group, and $\So=(\R\times
\hat{\Z})/\Z$.  $\So$ is a
`solenoid', is a compact, connected, abelian and uniquely divisible group.

\begin{proposition} With the above notation, 
\begin{enumerate}
\item   If $K$ is a number field, then $\phi$ is
surjective. 

\item $N_K$ is the connected component of identity in $W_K$, and there
is an isomorphism, 
\[\quad  N_K\simeq \Projlim_{\{L: K\subset
L\subset \bar{K}\}} D_L.  \]

\end{enumerate}
\end{proposition}

\section{Moore cohomology groups}
We recall here the cohomology theory using measurable cochains
developed by Moore in a series of papers (\cite{M1, M2, M3, M4},
especially  \cite{M3}). 

\begin{definition} A second countable topological group is said to be
{\em polish} if its topology admits a separable, complete metric. 
\end{definition}

Let $G$ be a separable, locally compact group. A $G$-module is a
topological abelian group with a continuous action $G\times A\to A$.
Denote by ${\mathcal P}(G)$  the category of commutative polish $G$-modules. 
For $A\in {\mathcal P}(G)$, let $I(A)=I(G, A)$ be the $G$-module obtained by
considering the `regular representation of $G$ with coefficients in
$A$': the underlying space of $I(A)$ is the collection of equivalence
classes of measurable functions from $G$ to $A$, where $G$ is equipped
with a Haar measure. $I(A)$ can 
be suitably topologized and the resulting topological module is again an
element of ${\mathcal P}(G)$. 

The existence theorem for a cohomology theory, is that there exists a
collection of covariant functors $H^r(G, .), ~r\geq 0, ~r\in Z$ from
${\mathcal P}(G)$ to the category of abelian groups satisfying the
following: 

\begin{description}
\item[C1] for every short exact sequence 
\[ 0\to A'\to A\to A''\to 0,\]
of modules in ${\mathcal P}(G)$, there is a long exact sequence of
cohomology groups, 
\[ \cdots \to H^{r-1}(G, A'')\to H^r(G, A')\to H^r(G, A)\to H^r(G,
A'')\to \cdots\]

\item[C2] $H^0(G, A)=A^G$, the space of $G$-invariants.

\item[C3] $H^r(G, I(A))=0$,  for all $r\geq 1$ and for all $A\in {\mathcal
P}(G)$. 

\end{description}

If $H_1$ and $H_2$ are two such cohomological functors
satisfying the above properties (1)-(3), it can be shown that 
there is a unique
isomorphism of functors from $H_1$ to $H_2$ compatible with the given
isomorphism in dimension $0$. 

For the existence part, let $C^n(G,A)$ be the set of all Borel
functions from $G^n$ to $A$. Denote by $\underline{C}^n(G,A)$ (also
equal to $I(G^n, A)$)
the set obtained from $C^n(G, A)$ by identifying $n$-cochains
which agree almost everywhere. Equipped with the usual coboundary
operator, the cohomology groups of these two complexes give raise to
two cohomological functors $H^n(G, A)$ and $\underline{H}^n(G, A)$,
satisfying C1-C3, and thus are isomorphic. With this isomorphism, the
cohomology groups $H^n(G, A)$ can be equipped with a topological
structure, though not necessarily Hausdorff. If further the topology
on $H^n(G, A)$  is Hausdorff, then it can be seen that $H^n(G, A)$ is
a polish group.

\begin{remark} Wigner \cite{Wi} has an alternate approach to these
cohomology groups, as Ext functors in a quasi-abelian category in the
sense of Yoneda. 
\end{remark} 

The low dimensional cohomology groups admit a concrete
interpretation. Let $Z^1(G, A)$ be the collection of continuous
`crossed homomorphisms' $c:G\to A$ satisfying $c(gh)=hc(g) +c(h)$ for
all $g, h\in G$. Denote by $B^1(G, A)$ the subset of $Z^1(G, A)$
consisting of the $1$-coboundaries of the form $c(g)=ga-a$ for some
$a\in A$. Then we have 
\[ H^1(G, A)\simeq Z^1(G, A)/B^1(G, A).\]
In particular, if the $G$-action is trivial on $A$, then 
$H^1(G,A)\simeq {\rm Hom}(G, A)$, the group of continuous 
homomorphisms from $G$ to $A$.  
In particular, the topology on $H^1(G, A)$ is the
topology give by the uniform convergence on compact subsets, and is
thus Hausdorff if $G$ acts trivially on $A$. 

For non-abelian coefficients $A$, the space   $H^1(G, A)$ can be
defined in a similar manner (see \cite[Annexe, Chapitre VII]{Se} and \cite[page 26]{M3}). Let 
\[ 1\to Z\to M\to M'\to 1,\]
be an exact sequence of locally compact topological groups, with $Z$
contained in the center of $M$. In exactly the same manner as in the
classical case, the following  sequence can be shown to be exact in
the appropriate sense:
\begin{equation}\label{nonabelianseq}
\cdots \to H^{1}(G, Z)\to H^1(G, M)\to H^1(G, M')\to H^2(G,
Z).
\end{equation}

One of the principal motivations for considering the cohomology theory
based on measurable cochains, is the isomorphism 
\[H^2(G, A)\simeq {\rm Ext}(G, A),\]
where ${\rm Ext}(G, A)$ is the set of equivalence classes of
extensions of topological groups 
\begin{equation}\label{extension}
 1\to A\to H\to G\to 1,
\end{equation}
equipped with the Baer product. Here $G$ is a separable locally
compact group and $A$ is polish. This follows essentially from a
theorem of Dixmier guaranteeing the 
existence of a Borel measurable cross-section $s:G\to H$ in the
situation of Equation \ref{extension}.  In
particular, we obtain 
\begin{proposition} Let $G$ be a separable, locally compact group and 
$A\in {\mathcal
P}(G)$. 
If $H^2(G, A)$ is trivial, then any exact sequence of topological
groups
\[ 1\to A\to H\to G\to 1\]
splits. 
\end{proposition}

There are many comparison theorems of the measurable cohomology
theory with continuous cohomology theories. If $G$ is discrete, then
the measurable, continuous and abstract cohomology theories are
isomorphic. It was shown by Wigner
\cite{Wi},  that if $G$ is a profinite group and $A$ is a discrete
$G$-module, then there is a natural isomorphism for all $i\geq 0$, 
\[ H^i_{gc}(G, A)\simeq H^i(G,A).\]
Consequently, we have
\begin{lemma}\label{wigner} 
Let $G$ be a profinite group acting algebraically on a torus $T$ defined
over $\C$. Then for $i\geq 1$,  
\[ H^i(G, T(\C))\simeq H^i_{gc}(G, T(\C)).\]
In particular, $ H^2(G, T(\C))\simeq H^2_{gc}(G, T(\C))$. 
\end{lemma}
\begin{proof}
Let  $X_*(T)$ denote the group of co-characters of $T$. Consider the
short exact sequence,
\begin{equation}\label{cocharacterseq}
1\to X_*(T)\to X_*(T)\otimes \C\to T(\C)\to 1.
\end{equation}
From the long exact sequences associated to the two cohomology
theories, and the fact that $H^1(G, V)$ vanishes for any vector space
$V$ over $\C$ in both the cohomology theories, we obtain for $i\geq 1$,
\[ H^i_{gc}(G, T(\C))\simeq H^{i+1}_{gc}(G, X_*(T)),\]
and similarly for the measurable cohomology groups. The lemma then
follows by Wigner's comparison theorem.
\end{proof}

The cohomology groups constructed by Moore have the added advantage of
the presence of a Hochshchild-Serre type spectral sequence. Let $G$ be
a locally compact group, $H$ be a closed normal subgroup of $G$ and
 $A \in {\mathcal P}(G)$.
\begin{proposition}[\cite{M3}[Theorem 9, page 29] \label{spectralsequence}
Assume further that $A$ is locally compact. 
There is a spectral sequence $E_r^{p, q}$ converging to $H^*(G,
A)$. If further $H^q(H, A)$ is Hausdorff, then we have for any $p$, 
\[ E_2^{p, q}= H^p(G/H, H^q(H, A)).\]
In particular, for $q=1 ~\text{or} ~0$, we have 
$ E_2^{p, q}= H^p(G/H, H^q(H, A)).$ 
\end{proposition}
The presence of the spectral sequence, together with the explicit form
of the $E_2^{p,q}$ terms for  $q=1 ~\text{or} ~0$ implies that the
following sequence is exact:
\[
\begin{split}\label{inflation-restriction}
1\to  H^1(G/H, A^H) & \xrightarrow{\text{inf}}  H^1(G,
 A)\xrightarrow{\text{res}}  H^0(G/H, H^1(H, A))\to \\
 &\to  H^2(G/H, A^H)\xrightarrow{\text{inf}}  H^2(G, A).
\end{split}
\]

We require the existence of the spectral sequence in the form of the 
following corollary:
\begin{corollary} In the notation of Proposition \ref{spectralsequence}, assume
further the following:
\begin{itemize}
\item the differential from $E_2^{1,1}\to E_2^{3,0}$ is injective. 
\item $H^2(H, A)=0$. 
\item the inflation map given by the spectral sequence 
\[ H^2(G/H,H^0(H, A)) \to H^2(G, A)\]
is zero. 
\end{itemize} 
Then $H^2(G, A)=(0)$. 
\end{corollary}

\section{Proof of the Main theorem for $\Z$-cases} 
We first consider the $\Z$-cases and compare
the cohomology of a discrete $G_K$-module with the cohomology of the
corresponding $W_K$-module:
\begin{proposition} Let $K$ be a non-archimedean local field or a
global field of positive characteristic ($\Z$-cases). Let $A$ be a
discrete $G_K$-module, and consider $A$ as a $W_K$-module via the
morphism $\phi$. Then for $i\geq 2$, the restriction map 
\[ \phi^*: H^i(G_K, A)\to H^i(W_K, A),\]
is an isomorphism.
\end{proposition}
\begin{proof} We have the following exact sequences of topological
groups:
\[ 1\to N_K\to G_K\to \hat{\Z}\to 1, \]
\[ 1\to N_K\to W_K\to {\Z}\to 1. \]
For a profinite group $G$, and a discrete $G$-module, the measurable
and continuous cohomology groups are isomorphic \cite{Wi}. The
cohomology group
\[H^q(N_K, A)=\Injlim_{U} H^q(N_K/U, A^U),\]
is a direct limit of the cohomology groups of the finite groups
$N_K/U$ as $U$ runs over the family of open subgroups of $N_K$. Thus
for $q\geq 1$, the cohomology groups $H^q(N_K, A)$ are torsion and discrete. 

We can
apply the spectral sequence, and compute it using the $E_2$ term of
the spectral sequence. The measurable cohomology groups of $\Z$ are
isomorphic to the abstract cohomology groups and thus  vanish in
degrees greater than one. For $\hat{\Z}$ (see \cite[Chapter  XIII,
Section 1]{Se}), 
the cohomology groups vanish
for degrees greater than one. Hence the Hochschild-Serre spectral sequence computing
the higher degree cohomology of $G= G_K$  or $W_K$, with the normal subgroup
being $N=N_K$  degenerates
at the $E_2$-stage, and for $n\geq 2$ we have the
restriction-inflation exact sequences,
\[ 1\to H^1(G/N, H^{n-1}(N, A))\to  H^n(G, A)\to H^0(G/N, H^n(N, A))\to 0. \]
\[\text{Now}\quad  H^1(\Z, A)= A/(F-1)A,\]
where  $F$ is the image of $1\in \Z$ in $\text{Aut}(A)$. 
Using the computation of the
cohomology of finite cyclic groups, we obtain for any discrete,
$\hat{\Z}$-module $A$,
\[H^1(\hat{\Z}, A)\simeq A'/(F-1)A,\]
where $A'=\{a\in A\mid (1+F+\dots F^n)a=0, \text{for some $n$}\}$. 
 Now it can be seen
that the torsion part
$A_t\subset A'$ \cite{Se}. Hence we have that if $A$ is also torsion
as an abelian group, then 
\[ H^1(\hat{\Z}, A)\simeq H^1(\Z, A).\] 
Since $H^i(N_K, A)$ is torsion for $i\geq 1$, we obtain the proposition.
\end{proof}

\begin{corollary} Let $G_K$ act algebraically on a torus $T$ defined
over $\C$,  and let $W_K$ act on $T$ via the induced action. Then the 
cohomology groups 
$H^i(W_K, T(\C))$ vanish for $i\geq 2$. In particular, 
Theorem \ref{maintheorem} is 
true in the $\Z$-cases.
\end{corollary}
\begin{proof} Consider the exact sequence given by Equation
(\ref{cocharacterseq}).  To see
that the cohomology groups $H^i(W_K, X_*(T)\otimes \C)$ vanish for
$i\geq 2$, we apply the Hochschild-Serre spectral sequence for the
Weil group $W_K$ and it's subgroup $N_K$. The only $E_2^{p,q}=
H^p(\Z, H^q(N_K, X_*(T)\otimes \C))$ terms that
are possibly non-vanishing are when $q=0$ and $p=0~\text{or}
~1$. Hence for $i\geq 2$, the groups $H^i(W_K, X_*(T)\otimes \C)$
vanish,  and we have $H^i(W_K, T(\C))\simeq H^{i+1}(W_K, X_*(T))$. 
By the above proposition and Lemma \ref{wigner}, we obtain for $i\geq 2$, 
\[H^i(W_K, T(\C))\simeq H^i(G_K, T(\C))\simeq H^i_{gc}(G_K, T(\C)).\]
The corollary then follows from the theorem of Poitou and  Tate (see
\cite{NSW}), 
asserting that the groups $G_K$ have strict (Galois) cohomological dimension 
$2$.
\end{proof}

\section{Proof of the Main theorem for archimedean local fields} 
 We first prove the following key proposition,
which is of use over  number fields too. 
\begin{proposition}\label{keyprop}
Let $H$ be a connected, locally compact abelian group, and  $S^1$
be  the circle group equipped with the trivial action of $H$. Then 
\[H^2(H, {{S}^1})=(0).\]
\end{proposition}
\begin{proof} The elements of $H^2(H, {{S}^1})$ parametrize equivalence
classes of extensions, 
\[ 1\to {{S}^1}\to G\to H\to 1,\] 
where ${{S}^1}$ is a closed normal subgroup of $G$, and $H$ is isomorphic
to the quotient $G/{{S}^1}$. We now use a construction due to Hughes
\cite{Hu}. Given a class $\alpha\in H^2(H, {{S}^1})$, 
and an extension  corresponding to $\alpha$ as above, choose a
measurable cross-section $\sigma: H\to G$, and define for $x, y\in H$
\[f_{\sigma}(x,y)=\sigma(x)\sigma(y)\sigma(x)^{-1}\sigma(y)^{-1}.\]
$f_{\sigma}(x,y)$ takes values in $A$, and it can be checked that
$f_{\sigma}(x,y)$ is a bilinear, measurable function from $H\times
H\to {{\mathbb S}^1}$. Now by Banach's theorem \cite{N}, 
any measurable homomorphism between
two locally compact groups is continuous, and so we obtain a
measurable homomorphism from $H\to {\rm Hom}(H, {{S}^1})$, where ${\rm
Hom}(H, {{S}^1})$ is the space of continuous, unitary homomorphisms of $H$,
and hence is topologically isomorphic to the dual group $\hat{H}$ of
$H$. But since $H$ is compact,  abelian, $\hat{H}$ is discrete. By
Banach's theorem, the homomorphism  $H\to \hat{H}$ is continuous.
Since $H$ is connected, this map has to be the trivial homomorphism.  

Hence we have that any extension of $H$ by ${S^1}$ is an abelian, locally
compact group. But the identity morphism of ${S^1}$, can be extended to a
morphism from $G\to { S^1}$, and hence any extension splits, and the
cohomology group is trivial.
\end{proof}

\subsection{Proof for $W_{\C}$}
\begin{corollary} The main theorem is true for $W_{\C}$. 
\end{corollary}
\begin{proof} This follows from the facts that for any compact group
$G$ and a topological vector space $V$, the higher cohomology groups
$H^i(G, V)$ vanish for $i\geq 1$, and the fact that $H^2(\R, \R)$
vanishes. 
\end{proof}

\subsection{Proof of Theorem \ref{maintheorem} for $W_{\R}$}
We have to show  for a $G(\C/\R)$-tori $T$, that $H^2(W_{\R},
T(\C))=(0)$. Let $\C^*_0, ~S^1_0, ~\Z_0$ denote respectively $\C^*,
~S^1, ~\Z$
with the trivial action of
$\Z/2\Z$, and let  $\C^*_1, ~S^1_1, ~\Z_1$ denote respectively $\C^*,
~S^1, ~\Z$,  where  $-1\in  \Z/2\Z$ acts
 by the non-trivial automorphism 
$z\to z^{-1}$. From the structure theory of torsion-free indecomposable
$\Z[\Z/2\Z]$-modules (\cite[Section 74, Chapter XI]{CR}), 
it follows that any torus $T$ with a $\Z/2\Z$-action
has a $\Z/2\Z$-equivariant 
filtration 
\[T(\C)=T_0(\C)\supset T_1(\C)\supset T_2\dots \supset (1),\]
 such that the subquotients are either  the  `invariant' $\C^*_0$, 
or the `anti-invariant'
$\C^*_1$.  Thus we can assume that we are in one of these two cases. 

We apply the Hochschild-Serre spectral sequence to the exact sequence,
\[ 1\to \C^*\to W_{\R}\to \Z/2\Z\to 1.\]
By Proposition \ref{keyprop}, we have that $E_2^{0,2}=H^0(\Z/2\Z,
H^2(\C^*, T(\C))=0$. 
From  the existence of the spectral sequence,  it follows that we have
to show the following for $j=0,1$: 
\begin{description}
\item[(SS1)] The differential 
$$d^0: H^0(\Z/2\Z, H^1(\C^*, \C^*_j)) \to H^2(\Z/2\Z, H^0(\C^*,
  \C^*_j))$$ 
is surjective.
\item[(SS2)] The differential $$d^1: H^1(\Z/2\Z, H^1(\C^*, \C^*_j)) \to 
 H^3(\Z/2\Z, H^0(\C^*, \C^*_j))$$ is injective.
\end{description}
The differential $E_2^{m-1,1}\to E_2^{m+1, 0}$ has been calculated by
Hochschild and Serre \cite{HS}[Theorem 4], in the situation when the normal
subgroup acts trivially on the coefficients.  The negative of the 
differential is
given by taking  cup-product with the class 
$a\in H^2(\Z/2\Z, \C^*_1)$, given  as in Equation
(\ref{extensionclass}) and  
defined  by the extension defining
$W_{\R}$. Since $\C^*$ acts trivially on $T(\C)$, we can use this
description of the differential, and further $ H^1(\C^*,
\C^*_j)=\text{Hom}(\C^*, \C^*_j)$. 

We claim the following: 

{\em Claim.} For any integer $i$, the cup-product 
\[\cup ~a:  \hat{H}^i(\Z/2\Z, \text{Hom}(S^1, S^1_j)) \to  
 \hat{H}^{i+2}(\Z/2\Z,  S^1_j)\]
is an isomorphism, where $\hat{H}^i$ stands for the Tate cohomology
groups. 

Write $\C^*=S^1\times \R$. The cohomology groups $H^i(\R, A)$ for any
coefficient group $A$ are uniquely divisible, and hence the Tate
cohomology groups $\hat{H}^j(\Z/2\Z, H^i(\R, A))$  vanish  for any
integer $j$. Hence the above claim implies the corresponding statement
of the claim with the
group $S^1$ replaced by $\C^*$, and the coefficients $S^1_j$ by
$\C^*_j$ for $j=0,1$.  
Since for any finite group $G$ and coefficient module $A$, we
have a surjection from $H^0(G,A)\to \hat{H}^0(G, A)$, it is clear that
the claim implies the above desired properties (SS1) and (SS2), 
of the differential in
the spectral sequence, and hence the theorem.

The proof of the claim is an explicit,  case by case calculation with
cup-products. By the definition of the
extension defining $W_{\R}$, we have that $-1\in \Z/2\Z$ acts by
$z\to z^{-1}$.  Since $\Z/2\Z$ acts on $S^1$ by $z\mapsto
z^{-1}$, we have the following isomorphisms as $\Z/2\Z$-modules:
\[ \text{Hom}(S^1, S^1_0)\simeq \Z_1 \quad \text{and}\quad 
\text{Hom}(S^1, S^1_1)\simeq \Z_0.\]

Let $\chi:\Z/2\Z \to \Q/\Z$ be the non-trivial
homomorphism, and let $\delta\chi\in H^2(\Z/2\Z,\Z)$ be the image with
respect to the connecting homomorphism of the exact seqeuence $0\to
\Z\to \Q\to\Q/\Z\to 0$. It follows from Herbrand's theorem that the
cup-product 
\[ \cup~\delta\chi: \hat{H}^i(\Z/2\Z,A)\to\hat{H}^{i+2}(\Z/2\Z,A)\]  is
an isomorphism for all values of $i$. $\delta\chi$ is given by the
formula:
\begin{equation*}
\delta\chi(\tau_1, \tau_2)=
\begin{cases} 1 & \text{if either $\tau_1$ or
$\tau_2$ is $1$},\\ $-1$ & \text{if $\tau_1=\tau_2=-1$}.
\end{cases}
\end{equation*}

The group $\Z/2\Z$ considered as a $\Z/2\Z$-module, 
occurs naturally as a quotient of $\Z_j$ and as a
subgroup of $S^1_j$ for $j=0,1$.  
Comparing the formulas for the cohomology class $a$ thought of as an
element in $ H^2(\Z/2\Z, S^1_1)$ and
$\delta\chi$, we obtain a commutative diagram, 
\[
\begin{CD} 
 \hat{H}^{i}(\Z/2\Z,  \Z_{j'})  @>>>  \hat{H}^{i}(\Z/2\Z,  \Z/2\Z) \\
@VV{\cup~a}V  @VV{\cup~\delta\chi}V\\
 \hat{H}^{i+2}(\Z/2\Z,  S^1_j) @<<<  \hat{H}^{i+2}(\Z/2\Z,  \Z/2\Z)
\end{CD}
\]
 where $j'\neq j$ and
$j',j\in \{0,1\}$. In order to prove the claim,  
we have to show that the left hand vertical map is
an isomorphism. From the formulas giving the Tate cohomology groups of
cyclic groups as in \cite[page 108]{CF}, we obtain natural isomorphisms
\begin{alignat*}{2}
 \hat{H}^{0}(\Z/2\Z,  \Z_{1}) & \simeq  \hat{H}^{0}(\Z/2\Z,
 \Z/2\Z) & \simeq  \hat{H}^{0}(\Z/2\Z,  S^1_{0})\\
 \hat{H}^{1}(\Z/2\Z,  \Z_{0}) & \simeq  \hat{H}^{1}(\Z/2\Z,
 \Z/2\Z) & \simeq  \hat{H}^{1}(\Z/2\Z,  S^1_{1}). 
\end{alignat*}{2}
In the other cases, the cohomology
groups $ \hat{H}^{i}(\Z/2\Z,  \Z_{j'})$ and $ \hat{H}^{i+2}(\Z/2\Z,
S^1_j)$ both vanish simultaneously.  The claim and Theorem
\ref{maintheorem} for $W_{\R}$ follows now from the commutativity  of
the above diagram and the fact that the cohomology of cyclic groups
are periodic of period two.

\begin{remark}\label{archremark}
 It is clear that the properties (SS1) and (SS2) hold
for any $T(\C)$ in place of $\C^*_0$ and $\C^*_j$, in view of the
structure theorem of torsion-free indecomposable
$\Z[\Z/2\Z]$-modules.
\end{remark}

\section{Proof of Main Theorem  for number fields}

We apply the Hochschild-Serre spectral sequence to the exact sequence,
\[ 1\to N_K\to W_{K}\to G_K\to 1.\]
By Proposition \ref{keyprop}, since the solenoid is a connected
abelian group, we have for any natural numbers $a, ~b$, 
\[H^2(\So^a\times (S^1)^b,  T(\C))=(0),\]
for the trivial action on $T(\C)$. Now $H^2(\R, \R)$ is zero, and
since the second cohomology commutes with projective limits of compact
groups \cite{M1}[Theorem 2.3], we conclude  that 
\[ E_2^{0,2}=H^0(G_K, H^2(N_K, T(\C))=(0).\]

In order to prove that the $E_2^{2,0}$ does not contribute to the
cohomology of the Weil group, 
we recall now a lemma proved by Langlands \cite{Lan2}[Lemma 4] and
by Labesse \cite{Lab}[Proposition 5.5(b)] (and providing for us the
analogue of Tate's vanishing theorem):  
\begin{lemma} \label{langlands}
Let $L$ be a finite extension of $K$ over which $T$ splits. The
inflation map,
\[ H^2(G_{L/K}, T(\C))\to H^2_{ct}(W_{L/K}, T(\C))\]
is zero, where the $H^*_{ct}$ denotes the continuous cohomology groups
of the arguments. 
\end{lemma}

We have the commutative diagram, 
\[
\begin{CD} 
{\Injlim_{\{M: L\subset M\subset \bar{K}\}} H^2(G_{M/K}, T(\C))}
 @>{\text{inf}}>>  H^2(G_{K}, T(\C)) \\
@VV{\text{inf}}V  @VV{\text{inf}}V\\
{\Injlim_{\{M: L\subset M\subset \bar{K}\}} H^2(W_{M/K}, T(\C))} 
@>{\text{inf}}>>  H^2(W_{K}, T(\C))
\end{CD}
\]
By Lemma \ref{wigner}, the top horizontal map is an
isomorphism. Since the left  hand side vertical map is zero by Lemma
\ref{langlands}, 
we obtain that the
inflation map 
\[H^2(G_{K}, T(\C))\to H^2(W_{K}, T(\C))\] 
is zero. 

Finally we are left with showing the differential $d_1: E_2^{1,1}\to
E_2^{3,0}$ is injective. We have
\begin{equation*}
\begin{split}
E_2^{1,1}& = H^1(G_K, \text{Hom}(N_K, T(\C)))\\
 & = \Injlim_{\{L:K\subset L\subset \tilde{K}\}} H^1(G_{L/K},
 \text{Hom}(D_L, T(\C))) 
\end{split}
\end{equation*}
From the structure theorem for $D_L$ given by Equation (\ref{DL}), we
have 
\[\text{Hom}(D_L, T(\C^*))\simeq  \text{Hom}(\R\times
\So^{r_1(L)+r_2(L)-1}, T(\C))\oplus  \text{Hom}({(S^1)}^{r_2(L)},
T(\C)).\]
Since the first group on the right hand side of the above equation is
uniquely divisible and $G_{L/K}$ is a finite group, we have 
\[
H^1(G_{L/K}, \text{Hom}(\R\times
\So^{r_1(L)+r_2(L)-1}, T(\C))=(0). \]
We need to analyze the structure of $\text{Hom}({(S^1)}^{r_2(L)},
T(\C^*))$ as a $G_K$-module. By passing to a larger extension we can
assume that $L$ has no real embeddings. Let $v$ be an archimedean
place of $K$. For each archimedean prime $v$ of $K$, fix a prime $w$
of $L$ dividing $v$. Let $G=G_{L/K}$ and let 
$G_v$ be the decomposition group of $G$ at $w$. 
Denote by $S_{\infty}(K)$ the set of archimedean places of $K$. 
It can be seen that, 
\[
\text{Hom}({(S^1)}^{r_2(L)},T(\C)) \simeq \oplus_{v\in
S_{\infty}(K)}\text{Ind}_{G_v}^G\text{Hom} (S^1, T(\C)).  \]
 Hence by Shapiro's lemma,
\[ H^1(G_{L/K},\text{Hom}({(S^1)}^{r_2(L)},
T(\C^*))\simeq \oplus_{v\in
S_{\infty, r}(K)}H^1(G_v, \text{Hom}(S^1, T(\C))\]
where $S_{\infty, r}(K)$ denotes the set of real places of $K$. In
each of the summands we have that complex conjugation $\sigma_v\in G_v$ acts
non-trivially on $S^1$. 
We thus obtain the following commutative diagram:
\[
\begin{CD} 
 H^1(G_{{L}/K}, H^1(D_L, T(\C))) @>{d_1}>>  H^3(G_K, T(\C))\\
@VV{\simeq}V  @VV{\text{res}}V\\
\oplus_{v\in S_{\infty, r}(K)}H^1(G_v, \text{Hom}(S^1, T(\C))  @>>> 
\oplus_{v\in S_{\infty, r}(K)}H^3(G_v, T(\C))
\end{CD}
\]
The nontrivial
element of the group $G_v$ for a real place $v$ of $K$,   acts on
$T(\C)$ by 
$\sigma_v$, and acts on $S^1$ by the inverse automorphism.
To show that the differential $d_1$ is injective, 
it is enough to show that the bottom arrow is injective.  But for each
real place $v$, this map is
the corresponding differential (see Property (SS2) and also Remark
\ref{archremark}), that came up in the computation of
$H^2(W_{\R}, T(\C))$.  Hence by
the computations done in the case of reals, we obtain that this map is
injective, and it follows that $E_3^{11}=(0)$. This proves Theorem
\ref{maintheorem} in the case of number fields.

\begin{remark} Over number fields, 
 it is not true in general that $H^2_{gc}(G_K, T(\C))$
vanishes. For
example, let  $K/\Q$ be  an imaginary quadratic extension of the
rationals and assume that complex conjugation acts by the non-trivial
automorphism of $\C^*$. It follows from 
the  results of Tate and Poitou, that in this situation 
$H^2(G_{\Q}, \C^*)\neq 0 $. 

It is tempting in the context of the results of this paper to
conjecture that the higher cohomology groups $H^i(W_K, T(\C))$ vanish
for $i\geq 3$. However this conjecture is false, as pointed out to
S. Lichtenbaum, even in the simplest situation when $K=\C$.  The Weil
group is essentially the circle group $S^1$ as far as cohomology is
concerned when $K=\C$.  By Theorem 4 of Wigner \cite{Wi}, the groups
$H^*(S^1,\Z)$ are isomorphic to the sheaf cohomology groups
$H^*(BS^1,\Z)$ of the classifying space $BS^1$ of $S^1$ with values in
the constant sheaf with structure group $\Z$ on $BS^1$. It is known
that the classifying space $BS^1$ of the circle group $S^1$ is
isomorphic to the infinite dimensional complex projective space
$\C\P^{\infty}=\Injlim_{n\in \N} \C\P^n$, and hence we obtain that
$H^i(S^1,\Z)$ is isomorphic to $\Z$ when $i$ is even and $0$
otherwise. Thus $H^i(W_{\C},\C^*)$ is isomorphic to $\Z$ when $i$ is
odd and vanishes otherwise.
\end{remark}

\section{Cohomology of Weil-Deligne groups}
Let $K$ be a local or global field, and let 
\[ W'_K=W_K\times SU_2\]
be the Weil-Deligne group, where $SU_2$ denotes the special unitary
group in two variables. We consider an action of $W'_K$ on a torus $T$
defined over $\C$ as before, via the natural projection to $G_K$ and
letting $G_K$ act through a finite quotient as automorphisms of
$T$. Consider $T(\C)$ equipped with the usual Euclidean topology as a $W_K$
module. We have, 
\begin{theorem} 
\[ H^2(W_K', T(\C))=(0).\]
\end{theorem}
\begin{proof} 
By the results of Moore (see \cite{M5} and \cite{M4}), it is known
that a topologically simply connected Lie group $G$ satisfies,
\[ H^1(G, S^1)=(0) \quad \text{and}\quad  H^2(G, S^1)=(0).\]
Since $SU_2$ is compact, the cohomology groups $H^i(SU_2, V)$ vanish
for $i\geq 1$ and for any real vector space $V$. Hence we obtain, 
\[ H^1(SU_2, T(\C))=(0) \quad \text{and}\quad  H^2(SU_2, T(\C))=(0).\]
Applying the Hochschild-Serre-Moore spectral sequence, with the normal
subgroup as $SU_2$, the spectral sequence collapses. The Main theorem 
then yields the above theorem. 
\end{proof}

\section{Application to lifting: Labesse's theorem and base change}
The principal motivation for Theorem \ref{maintheorem} is the
application to lifting admissible homomorphisms of Weil groups. 
 Let $G$ be a reductive group defined over $K$.  Denote by 
${^L}{{G}}^0$ the complex points of  the Langlands dual group 
defined over the complex numbers. The Galois group $G_K$ acts on
${^L}{{G}}^0$, and  let ${^L}{{G}}$ the semi-direct
product of ${^L}{{G}}^0$ by $W_K$, where $W_K$ acts via the natural
projection to $G_K$.    An admissible
homomorphism $\phi: W_K'\to {^L}{{G}}$ consists of the following
(for further details we refer to \cite{Bo} and \cite{BR}):
\begin{itemize} 
\item $\phi$ is continuous and composition with the projection to
$W_K$ induces the identity on $W_K$. 
\item for every $w\in W_K$, $\phi(w)$ is semisimple.
\end{itemize}
Let $A(W_K', {^L}{{G}}^0)$ denote the collection of admissible
homomorphisms. It can be seen that $A(W_K', {^L}{{G}}^0)\subset
H^1(W_K', {^L}{{G}}^0)$. 
Assume now that there is an
exact sequence of the form,
\begin{equation}\label{lgroupextension}
1\to S \to {^L}{\tilde{G}}^0\xrightarrow{\pi} {{^L}{G}^0} \to 1,
\end{equation}
where $G, \tilde{G}$ are as above, $S$ is a central torus in ${^L}{\tilde{G}}^0$, and the
morphism $\pi$ is compatible with the action of $G_K$. As a corollary
of Theorem \ref{maintheorem} and the exact sequence
Equation (\ref{nonabelianseq}),   we have
\begin{Labesse} \label{labesse}
With notation as above, the natural map
\[ A(W_K', {^L}{\tilde{G}}^0) \to A(W_K', {^L}{{G}}^0)\]  
is surjective.
\end{Labesse}
\begin{remark} It is also possible to see this directly from our Main
theorem as follows: choose a measurable cross-section $s:{^L}{{G}}^0
\to {^L}{\tilde{G}}^0$. Corresponding to an admissible map $\phi\in
A(W_K', {^L}{{G}}^0)$, we obtain the pullback of the above extension
given by Equation (\ref{lgroupextension}):
\begin{equation*}
1\to S \to E \to W_K' \to 1. 
\end{equation*}   
Hence our Main theorem implies the existence of an admissible lift. 
\end{remark} 
\begin{remark} Let ${\mathcal A}$ be a von Neumann algebra realizable
on a separable Hilbert space. It follows from Theorem \ref{gentate}
and \cite[Theorem 5, page 40]{M4}, that any strongly continuous
representation $\rho$ of $W_K'$ by inner automorphisms of ${\mathcal
A}$ is implementable, i.e., there exists a continuous unitary
representation $\pi: W_K'\to U({\mathcal A})$, where $U({\mathcal A})$
is the unitary group in ${\mathcal A}$ such that for any $a\in
{\mathcal A}$ and $w\in W_K'$, we have
$\rho(w)(a)=\pi(w)a\pi(w)^{-1}$. It would be interesting to know of
any arithmetical consequences of this result.  
\end{remark} 

\subsection{Characterising image of base change}
Following earlier work by Saito and Shintani, 
Langlands in \cite{Lan1}, established the existence 
of the base change lift  $BC_{L/K}$ sending cuspidal automorphic
representations of $GL_2(\mathbf{A}_k)$ to automorphic representations of 
$GL_2(\mathbf{A}_L)$,  provided $L/K$ is a cyclic extension of
prime degree.  Langlands further characterized the image and the
fibres of the base change map $BC_{L/K}$. 

The work of Langlands was extended to $GL_n,
~n\geq 3$ by Arthur and Clozel \cite{AC}. 
The descent property of automorphic representations established by 
 Langlands, Arthur and Clozel is that 
 if  $\Pi$ is a cuspidal automorphic representation of
$GL_n(\mathbf{A}_L)$, which is invariant with respect to the action of 
the Galois group $G(L/K)$ of a cyclic extension $L/K$ of prime degree,
 then $\Pi$ lies in the image of the base change map
$BC_{L/K}$ 
from automorphic representations on $GL_n(\mathbf{A}_k)$ to automorphic
representations on $GL_n(\mathbf{A}_L)$.
The descent property for invariant representations with respect to a
cyclic group of automorphisms,  was obtained  Lapid
and Rogawski  for $GL(2)$, and conditionally by them for $GL(n), ~n\geq 3$.

We assume the functoriality principle, and consider the conjectural
complex points of the Langlands group,  the  representations of which parametrize automorphic
representations, as a generalization of the Galois (or Weil) group.
We have by  Tate's theorem for the Galois group $G_K$, and by the main
theorem of this paper for the Weil-Deligne group, that $ H^2(W_K',\C^*)=1.$
On the functorial side, base change from $K$ to $L$, 
 amounts to restriction to the 
subgroup $W_L$ of an admissible parameter in $A(W_K', GL_n(\C))$, since
$(^LGL_n)^0\simeq GL_n(\C)$. 
Conversely if   $\rho$ is  a $G(L/K)$-invariant representation of
the Galois group $W_L$ into $GL_n(\C)$, then the question of
characterising the image of base change for $GL_n$ from $K$ to $L$, 
translates to the question whether $\rho$ can
be extended to a representation of $G_K$. If the commutant of the
image of $\rho$ consists only of scalars, i.e., if $\rho$ is assumed to be
irreducible, then it can be 
extended as  a projective representation of $W_K'$ to $PGL(n,\C)$. By
Theorem \ref{gentate}, this projective representation has a lift to $GL(n,
\C)$. Restricting back to $W_L$, it can be checked that the lift
differs from the original representation by an invariant character. 

In the automorphic context, 
let $L/K$ be a  solvable extension of global fields, and let $\Pi$ be
a unitary, cuspidal automorphic representation of  
$GL_n(\mathbf{A}_L)$ which is ${\rm
Gal}(L/K)$-invariant. It is shown in \cite{R},  
using the results of Lapid and Rogawski \cite{LR} (based 
upon the assumption that a relevant fundamental
lemma holds for $n\geq 3$, and known to be true  when $n=2$),  that 
there exists a $G(L/K)$-invariant Hecke character
$\psi$ of $K$, 
 and a cuspidal automorphic representation $\pi$ of  
$GL_n(\mathbf{A}_L)$ such that 
\[ BC_{L/K}(\pi) \simeq \Pi\otimes \psi.\]
Further $\psi$ is unique upto base change to $L$ of a Hecke
character of $k$. 

A difficult case to handle in the proof of this theorem, say for $n=2$, 
 was when $\Pi$ is an automorphically induced representation. 
In this situation, the
 automorphic representations are parametrised by admissible parameters
of the Weil group, and the lifting theorem helps in producing 
the required descent (or extension) of the admissible parameter. It
 can then be checked that the extended parameter gives raise to a
 cuspidal representation.  Although the results of \cite{R} do not use
 this reasoning to obtain the result in the automorphic context, it is
 this argument that led us to prove 
the results of this paper.

\end{document}